\documentstyle[12pt]{article}
\begin{document}
\setlength{\textheight}{574pt}
\setlength{\textwidth}{432pt}
\setlength{\oddsidemargin}{18pt}
\setlength{\topmargin}{14pt}
\setlength{\evensidemargin}{18pt}
\newtheorem{theorem}{Theorem}[section]
\newtheorem{lemma}{Lemma}[section]
\newtheorem{corollary}{Corollary}[section]
\newtheorem{remark}{Remark}[section]
\newtheorem{sublemma}{Sublemma}[section]
\newtheorem{definition}{Definition}[section]
\newtheorem{problem}{Problem}
\newtheorem{proposition}{Proposition}[section]
\title{{\bf EFFECTIVE BIRATIONALITY OF PLURICANONICAL SYSTEMS}}
\date{June 2000}
\author{Hajime TSUJI}
\maketitle
\begin{abstract}
By using the theory of AZD originated by the author, I prove that 
for every smooth projective $n$-fold $X$ of general type and every 
\[
m\geq \lceil\sum_{\ell =1}^{n}\sqrt[\ell]{2}\,\,\ell\rceil +1,
\]
$\mid mK_{X}\mid$ gives a birational rational map from $X$
 into a projective space,
unless it has a nontrivial (relative dimension is positive) rational fiber space structure whose general fiber is birational to a variety of relatively low degree 
in a projective space. MSC 32J25 
\end{abstract}
\tableofcontents
\section{Introduction}
Let $X$ be a smooth projective variety and let $K_{X}$ be the canonical 
bundle of $X$. 
$X$ is said to be of general type, if $K_{X}$ is big, i.e.,
\[
\limsup_{m\rightarrow\infty}m^{-\dim X}\dim H^{0}(X,{\cal O}_{X}(mK_{X}))
> 0
\]
holds.
The following problem is fundamental to study projective 
vareity of general type. \vspace{5mm}\\
{\bf Problem}
Let $X$ be a smooth projective variety of general type.
Find a positive integer $m_{0}$ such that for every $m\geq m_{0}$,
$\mid mK_{X}\mid$ gives a birational rational map from $X$
into a projective space. \vspace{5mm} \\
If $\dim X = 1$, it is well known that $\mid 3K_{X}\mid$ gives a 
projective embedding.
In the case of smooth projective surfaces of general type, 
E. Bombieri showed that $\mid 5K_{X}\mid$ gives a birational
rational map from $X$ into a projective space (\cite{b3}).
In the case of $\dim X \geq 3$, 
I have proved the following theorem.
\begin{theorem}(\cite{tu4})
There exists a positive integer $\nu_{n}$ which depends
only on $n$ such that for every smooth projective $n$-fold $X$
of general type defined over complex numbers, $\mid mK_{X}\mid$ gives a birational rational map
from $X$ into a projective space for every $m\geq \nu_{n}$.
\end{theorem}
Theorem 1.1 is an affirmative answer to the problem. 
But it seems to be very hard to give an {\bf effective estimate} 
of the number $\nu_{n}$ because the proof depends on 
the abstract facts of Hilbert scheme. 

The main purpose of this article is to give the
following {\bf weak effective answer} to  the problem. 
\begin{theorem}
For every smooth projective $n$-fold $X$ of general type, one of the
followings holds.
\begin{enumerate}
\item for every 
\[
m\geq \lceil\sum_{\ell =1}^{n}\sqrt[\ell]{2}\,\,\ell\rceil +1,
\]
$\mid mK_{X}\mid$ gives a birational rational map from $X$
 into a projective space,
\item $X$ is dominated by a family of  subvarieties of dimension $d (\geq 1)$
which are 
birational to subvarieties of degree less than or equal to 
$(\lceil\sum_{\ell=1}^{n}\sqrt[\ell]{2}\,\,\ell\rceil +1)^{d}$
in a projective space by some pluricanonical system $\mid \alpha K_{X}\mid$.
\end{enumerate}
\end{theorem}

\begin{theorem}
For every smooth projective $n$-fold $X$ of general type, one of the
followings holds.
\begin{enumerate}
\item for every 
\[
m\geq \lceil\sum_{\ell =1}^{n}\sqrt[\ell]{2}\,\,\ell\rceil +1,
\]
$\mid mK_{X}\mid$ gives a birational rational map from $X$
 into a projective space,
\item there exists a rational fibration 
\[
f : X -\cdots\rightarrow Y
\]
such that a general fiber $F$ of $f$ is positive dimensional and 
is birational to a subvariety of degree less than or equal to 
$(\lceil\sum_{\ell=1}^{n}\sqrt[\ell]{2}\,\,\ell\rceil +1)^{d^{2}}d^{d}$ in a projective space by some pluricanonical system $\mid\alpha K_{X}\mid$, 
where $d$ denotes the dimension of $F$. 
\end{enumerate}
\end{theorem}

The proofs  of Theorem 1.2 and 1.3 are technically 
much easier than that of Theorem 1.1 (\cite{tu4}). 
But they 
are effective and 
clarify the essential obstruction to obtain the 
birationality of the pluricanonical map $\Phi_{\mid mK_{X}\mid}$
with relatively small $m$.
As one see in the proof, if we need very large $m$ to embed 
$X$ birationally into a projective space by $\mid mK_{X}\mid$,
the image $\Phi_{\mid mK_{X}\mid}(X)$ is distorted in the sense that 
$X$ is small in the fiber direction and large in the horizontal direction 
with respect to a rational fiber space structure.   
Such a phenomenon was first observed by E. Bombieri in his paper 
\cite{b3}. 
Actually he  found the existence of a genus 2 fibration is an obstruction 
to the birationality of $\mid 2K_{X}\mid$ for some surfaces of 
general type (\cite[p. 173, Main Theorem (iv)]{b3}). 
Of course the results above will not be optimal and 
more abstract in comparison with the case of surfaces.

In the case of 3-folds of general type, there were several results\cite{luo,luo2,b-d} in 
this direction.
But these results depend on the plurigenus formula 
for canonical 3-folds of general type and 
moreover their estimates depend on the apriori bound of 
$\chi (X,{\cal O}_{X})$, hence it is even weaker than 
Theorem 1.1 in this respect and the bound is 
is so huge that they have only a theoretical interests. 

I hope the estimates in Theorem 1.2 and Theorem 1.3 
are acceptable at least for projective varieties of low dimension.  
And it seems to be more or less optimal in order of size,
 even in the case of arbitrary dimension. 
But I should  say that even in the case of 3-folds, the exceptional cases 
seem to be very hard to classify. 

As in \cite{tu4}, the main difficulty is the fact that $K_{X}$ is not ample in general.
To overcome this difficulty we use a special singular hemitian metric
on $K_{X}$ called AZD which was originated by the author (\cite{tu}).
By using AZD we can handle $K_{X}$ as if $K_{X}$ were nef and big.

\section{Preliminaries}
\subsection{Multiplier ideal sheaves}
In this section, we shall review the basic definitions and properties
of multiplier ideal sheaves.
\begin{definition}
Let $L$ be a line bundle on a complex manifold $M$.
A singular hermitian metric $h$ is given by
\[
h = e^{-\varphi}\cdot h_{0},
\]
where $h_{0}$ is a $C^{\infty}$-hermitian metric on $L$ and 
$\varphi\in L^{1}_{loc}(M)$ is an arbitrary function on $M$.
\end{definition}
The curvature current $\Theta_{h}$ of the singular hermitian line
bundle $(L,h)$ is defined by
\[
\Theta_{h} := \Theta_{h_{0}} + \sqrt{-1}\partial\bar{\partial}\varphi ,
\]
where $\partial\bar{\partial}$ is taken in the sense of a current.
The $L^{2}$-sheaf ${\cal L}^{2}(L,h)$ of the singular hermitian
line bundle $(L,h)$ is defined by
\[
{\cal L}^{2}(L,h) := \{ \sigma\in\Gamma (U,{\cal O}_{M}(L))\mid 
\, h(\sigma ,\sigma )\in L^{1}_{loc}(U)\} ,
\]
where $U$ runs opens subsets of $M$.
In this case there exists an ideal sheaf ${\cal I}(h)$ such that
\[
{\cal L}^{2}(L,h) = {\cal O}_{M}(L)\otimes {\cal I}(h)
\]
holds.  We call ${\cal I}(h)$ {\bf the multiplier ideal sheaf of $(L,h)$}.
If we write $h$ as 
\[
h = e^{-\varphi}\cdot h_{0},
\]
where $h_{0}$ is a $C^{\infty}$ hermitian metric on $L$ and 
$\varphi\in L^{1}_{loc}(M)$ is the weight function, we see that
\[
{\cal I}(h) := {\cal L}^{2}({\cal O}_{M},e^{-\varphi})
\]
holds. 
We also denote ${\cal L}^{2}({\cal O}_{M},e^{-\varphi})$ 
by ${\cal I}(\varphi )$. 
Let $(L,h)$ be a singular hermitian line bundle on a smooth projective variety 
$X$ such that
\[
\Theta_{h} \geq - \omega
\] 
holds for some $C^{\infty}$ K\"{a}hler form $\omega$ on $X$.
Then by \cite[p. 561]{n}, we see that ${\cal I}(h)$ is a coherent sheaf of 
${\cal O}_{X}$-ideal. 

Similarly we obtain the sheaf 
\[
{\cal I}_{\infty}(h):= {\cal L}^{\infty}({\cal O}_{M},e^{-\varphi})
\]
and call it {\bf the $L^{\infty}$-multiplier ideal sheaf} of 
$(L,h)$. 
We have the following vanishing theorem.

\begin{theorem}(Nadel's vanishing theorem \cite[p.561]{n})
Let $(L,h)$ be a singular hermitian line bundle on a compact K\"{a}hler
manifold $M$ and let $\omega$ be a K\"{a}hler form on $M$.
Suppose that $\Theta_{h}$ is strictly positive, i.e., there exists
a positive constant $\varepsilon$ such that
\[
\Theta_{h} \geq \varepsilon\omega
\]
holds.
Then ${\cal I}(h)$ is a coherent sheaf of ${\cal O}_{M}$-ideal 
and for every $q\geq 1$
\[
H^{q}(M,{\cal O}_{M}(K_{M}+L)\otimes{\cal I}(h)) = 0
\]
holds.
\end{theorem}

\subsection{Analytic Zariski decomposition}

To study a big line bundle we introduce the notion of analytic Zariski
decompositions.
By using analytic Zariski decompositions, we can handle big line bundles
like a nef and big line bundles.
\begin{definition}
Let $M$ be a compact complex manifold and let $L$ be a line bundle
on $M$.  A singular hermitian metric $h$ on $L$ is said to be 
an analytic Zariski decomposition, if the followings hold.
\begin{enumerate}
\item $\Theta_{h}$ is a closed positive current,
\item for every $m\geq 0$, the natural inclusion
\[
H^{0}(M,{\cal O}_{M}(mL)\otimes{\cal I}(h^{m}))\rightarrow
H^{0}(M,{\cal O}_{M}(mL))
\]
is isomorphim.
\end{enumerate}
\end{definition}
\begin{remark} If an AZD exists on a line bundle $L$ on a smooth projective
variety $M$, $L$ is pseudoeffective by the condition 1 above.
\end{remark}

\begin{theorem}(\cite{tu,tu2} see also \cite[Section 2.2]{tu4})
 Let $L$ be a big line  bundle on a smooth projective variety
$M$.  Then $L$ has an AZD. 
\end{theorem}

\subsection{Volume of subvarieties}

To measure the positivity of big line bundles 
on a projective variety we shall introduce a volume of a projective variety
with respect to a line bundle.

\begin{definition} Let $L$ be a line bundle on a compact complex 
manifold $M$ of dimension $n$. 
We define the $L$-volume of $M$ by
\[
\mu (M,L) := n!\cdot\overline{\lim}_{m\rightarrow\infty}m^{-n}
\dim H^{0}(M,{\cal O}_{M}(mL)).
\]
\end{definition}
\begin{definition}(\cite{tu3})
Let $(L,h)$ be a singular hermitian line bundle on a smooth projective variety
$X$ such that $\Theta_{h}\geq 0$. 
Let $Y$ be a subvariety of $X$ of dimension $r$.
We define the volume $\mu (Y,L)$ of $Y$ with respect to 
$L$ by 
\[
\mu (Y,L) := r!\cdot\overline{\lim}_{m\rightarrow\infty}m^{-r}
\dim H^{0}(Y,{\cal O}_{Y}(mL)\otimes{\cal I}(h^{m})/tor),
\]
where $tor$ denotes the torsion part of the sheaf ${\cal O}_{Y}(mL)\otimes{\cal I}(h^{m})$. 
\end{definition}
\section{Stratification of varieties by multiplier ideal sheaves}
In this section we shall construct a stratification 
of a smooth projective variety $X$ of general type  by using an AZD $h$
of $K_{X}$.
We use  the ideas 
in \cite{a-s,t} to construct the stratification. 
But since $(K_{X},h)$ is not an ample line bundle, 
the argument is a little bit more involved.  
\subsection{Construction of a stratification} 
Let $X$ be a smooth projective $n$-fold of general type.
Let $h$ be an AZD of $K_{X}$. 
Let us denote $\mu (X,K_{X})$ by $\mu_{0}$.
We set 
\[
X^{\circ} = \{ x\in X\mid x\not{\in} \mbox{Bs}\mid mK_{X}\mid \mbox{and  
$\Phi_{\mid mK_{X}\mid}$ is a biholomorphism}
\]
\[
\hspace {50mm} \mbox{ on a neighbourhood of  $x$ 
for some $m\geq 1$}\} .
\]
Then $X^{\circ}$ is a nonempty Zariski open subset of $X$. 
\begin{lemma} Let $x,x^{\prime}$ be distinct points on $X^{\circ}$.  
We set 
\[
{\cal M}_{x,x^{\prime}} = {\cal M}_{x}\otimes{\cal M}_{x^{\prime}}
\]

Let $\varepsilon$ be a sufficiently small positive number.
Then 
\[
H^{0}(X,{\cal O}_{X}(mK_{X})\otimes{\cal M}_{x,x^{\prime}}^{\lceil\sqrt[n]{\mu_{0}}
(1-\varepsilon )\frac{m}{\sqrt[n]{2}}\rceil})\neq 0
\]
for every sufficiently large $m$, where ${\cal M}_{x},{\cal M}_{x^{\prime}}$ denote the
maximal ideal sheaf of the points $x,x^{\prime}$ respectively.
\end{lemma}
{\bf Proof of Lemma 3.1}.   
Let us consider the exact sequence:
\[
0\rightarrow H^{0}(X,{\cal O}_{X}(mK_{X})\otimes
{\cal M}_{x,x^{\prime}}^{\lceil\sqrt[n]{\mu_{0}}(1-\varepsilon )\frac{m}{\sqrt[n]{2}}\rceil})
\rightarrow H^{0}(X,{\cal O}_{X}(mK_{X}))\rightarrow
\]
\[
  H^{0}(X,{\cal O}_{X}
(mK_{X})/{\cal M}_{x,x^{\prime}}^{\lceil\sqrt[n]{\mu_{0}}(1-\varepsilon )\frac{m}{\sqrt[n]{2}}\rceil}).
\]
Since 
\[
n!\cdot\overline{\lim}_{m\rightarrow\infty}m^{-n}\dim H^{0}(X,{\cal O}_{X}(mK_{X})
/{\cal M}_{x,x^{\prime}}^{\lceil\sqrt[n]{\mu_{0}}(1-\varepsilon )\frac{m}{\sqrt[n]{2}}\rceil})
=
\mu_{0}(1-\varepsilon )^{n} < \mu_{0}
\]
hold, we see that Lemma 3.1 holds.  {\bf Q.E.D.}

\vspace{10mm}

Let us take a sufficiently large positive integer $m_{0}$ and let $\sigma$
be a general (nonzero) element of  
$H^{0}(X,{\cal O}_{X}(m_{0}K_{X})\otimes
{\cal M}_{x,x^{\prime}}^{\lceil\sqrt[n]{\mu_{0}}(1-\varepsilon )\frac{m_{0}}{\sqrt[n]{2}}\rceil})$.
We define a singular hermitian metric $h_{0}$ on $K_{X}$ by
\[
h_{0}(\tau ,\tau ) := \frac{\mid \tau\mid^{2}}{\mid \sigma\mid^{2/m_{0}}}.
\]
Then 
\[
\Theta_{h_{0}} = \frac{2\pi}{m_{0}}(\sigma )
\]
holds, where $(\sigma )$ denotes the closed positive current 
defined by the divisor $(\sigma )$.
Hence $\Theta_{h_{0}}$ is a closed positive current.
Let $\alpha$ be a positive number and let ${\cal I}(\alpha )$ denote
the multiplier ideal sheaf of $h_{0}^{\alpha}$, i.e.,
\[
{\cal I}(\alpha ) = 
{\cal L}^{2}({\cal O}_{X},(\frac{h_{0}}{h_{X}})^{\alpha}),
\]
where $h_{X}$ is an arbitrary $C^{\infty}$-hermitian metric on
$K_{X}$.
Let us define a positive number $\alpha_{0} (= \alpha_{0}(x,y))$ by
\[
\alpha_{0} := \inf\{\alpha > 0\mid ({\cal O}_{X}/{\cal I}(\alpha ))_{x}\neq 0\,\mbox{and}\, ({\cal O}_{X}/{\cal I}(\alpha))_{x^{\prime}}\neq 0\}.
\]
Since $(\sum_{i=1}^{n}\mid z_{i}\mid^{2})^{-n}$ is not locally integrable 
around $O\in \mbox{{\bf C}}^{n}$, by the construction of $h_{0}$, we see
that 
\[
\alpha_{0}\leq \frac{n\sqrt[n]{2}}{\sqrt[n]{\mu_{0}}(1-\varepsilon )}
\]
holds.
Then one of the following two cases occurs. \vspace{10mm} \\
{\bf Case} 1.1:  For every small positive number $\delta$, 
${\cal O}_{X}/{\cal I}(\alpha_{0}-\delta )$  has $0$-stalk 
at both $x$ and $x^{\prime}$. \\
{\bf Case} 1.2: For every small positive number $\delta$, 
${\cal O}_{X}/{\cal I}(\alpha_{0}-\delta )$  has nonzero-stalk 
at one of $x$ or $x^{\prime}$ say $x^{\prime}$. \vspace{10mm} \\

We first consider Case 1.1.
Let $\delta$ be a sufficiently small positive number and  
let $V_{1}$ be the germ of subscheme at $x$ defined by the ideal sheaf 
${\cal I}(\alpha_{0}+\delta )$.
By the coherence of ${\cal I}(\alpha ) (\alpha > 0)$, we see that 
if we take $\delta$ sufficiently small, then $V_{1}$ is independent
of $\delta$.  It is also easy to verify that $V_{1}$ is reduced if 
we take $\delta$ sufficiently small. 
In fact if we take a log resolution of 
$(X,\frac{\alpha_{0}}{m_{0}}(\sigma ))$, 
$V_{1}$ is the image of the divisor with discrepancy $-1$ 
(for example cf. \cite[p.207]{he}). 
Let $X_{1}$ be a subvariety of $X$ which defines a branch of $V_{1}$
at $x$. 
We consider the following two cases. \vspace{10mm} \\
{\bf Case} 2.1: $X_{1}$ passes through both $x$ and $x^{\prime}$, \\
{\bf Case} 2.2: Otherwise \vspace{10mm} \\

For the first we consider Case 2.1.
Suppose that $X_{1}$ is not isolated at $x$.  Let $n_{1}$ denote
the dimension of $X_{1}$.  Let us define the volume $\mu_{1}$ of $X_{1}$
with respect to $K_{X}$ by
\[
\mu_{1} := \mu (X_{1},K_{X}).
\]
Since $x\in X^{\circ}$, we see that $\mu_{1} > 0$ holds.

\begin{lemma} Let $\varepsilon$ be a sufficiently small positive number and let $x_{1},x_{2}$ be distinct regular points on $X_{1}\cap X^{\circ}$. 
Then for a sufficiently large $m >1$,
\[
H^{0}(X_{1},{\cal O}_{X_{1}}(mK_{X})\otimes{\cal I}(h^{m})\otimes
{\cal M}_{x_{1},x_{2}}^{\lceil\sqrt[n_{1}]{\mu_{1}}(1-\varepsilon )\frac{m}{\sqrt[n_{1}]{2}}\rceil})\neq 0
\]
holds.
\end{lemma}
The proof of Lemma 3.2 is identical as that of Lemma 3.1, since 
\[
{\cal I}(h^{m})_{x_{i}} = {\cal O}_{X,x_{i}} (i=1,2)
\]
hold for every $m$ by Proposition 2.1 and Lemma 2.1.

By Kodaira's lemma there is an effective {\bf Q}-divisor $E$ such
that $K_{X}- E$ is ample.
Let $\ell_{1}$ be a sufficiently large positive integer which will be specified later  such that
\[
L_{1} := \ell_{1} (K_{X}- E)
\]
is Cartier. 
\begin{lemma}
If we take $\ell_{1}$ sufficiently large, then 
\[
\phi_{m} : H^{0}(X,{\cal O}_{X}(mK_{X}+L_{1})\otimes{\cal I}(h^{m}))\rightarrow 
H^{0}(X_{1},{\cal O}_{X_{1}}(mK_{X}+L_{1} )\otimes{\cal I}(h^{m}))
\]
is surjective for every  $m\geq 0$.
\end{lemma}
{\bf Proof}.
Let us take a locally free resolution of the ideal sheaf ${\cal I}_{X_{1}}$
of $X_{1}$.
\[
0\leftarrow {\cal I}_{X_{1}}\leftarrow {\cal E}_{1}\leftarrow {\cal E}_{2}
\leftarrow \cdots \leftarrow {\cal E}_{k}\leftarrow 0.
\]
Then by the trivial extension of  
the case of vector bundles, if $\ell_{1}$ is sufficiently large, we see that
\[
H^{q}(X,{\cal O}_{X}(mK_{X}+L_{1})\otimes{\cal I}(h^{m})
\otimes{\cal E}_{j}) = 0
\]
holds for every $m\geq 1$, $q\geq 1$ and  $1\leq j\leq k$.
In fact if we take $\ell_{1}$ sufficiently large, we see that for every $j$, 
${\cal O}_{X}(L_{1} - K_{X})\otimes {\cal E}_{j}$ admits a $C^{\infty}$-hermitian metric $g_{j}$ such that
\[
\Theta_{g_{j}} \geq \mbox{Id}_{E_{j}}\otimes \omega
\]
holds, where $\omega$ is a K\"{a}hler form on $X$.
By \cite[Theorem 4.1.2 and Lemma 4.2.2]{ca} we have the desired vanishing. 
Hence we have:
\begin{sublemma}
\[
H^{1}(X,{\cal O}_{X}(mK_{X}+L_{1})\otimes{\cal I}(h^{m})\otimes{\cal E}_{j}) = 0
\]
holds for every $m\geq 0$ and $1\leq j\leq r$
\end{sublemma}
Let 
\[
p_{m} : X_{m}\longrightarrow X
\]
be a composition of successive blowing ups with smooth centers 
such that 
$p_{m}^{*}{\cal I}(h^{m})$ is locally free on $X_{m}$. 
\begin{sublemma}
\[
R^{p}p_{m*}({\cal O}_{X_{m}}(K_{X_{m}})\otimes {\cal I}(p_{m}^{*}h^{m})) = 0
\]
holds for every $p\geq 1$ and $m\geq 1$. 
\end{sublemma}
We note that by the definition of multiplier ideal sheaves 
\[
p_{m*}({\cal O}_{X_{m}}(K_{X_{m}})\otimes {\cal I}(p_{m}^{*}h^{m}))
= {\cal O}(K_{X})\otimes {\cal I}(h^{m})
\]
holds. 
Hence by Sublemma 3.1 and Sublemma 3.2 and the Leray 
spectral sequence, we see that 
\[
H^{q}(X_{m},{\cal O}_{X_{m}}(K_{X_{m}}+p_{m}^{*}(mK_{X}+L_{1} - K_{X}))
\otimes {\cal I}(p_{m}^{*}h^{m})\otimes p_{m}^{*}{\cal E}_{j}) = 0
\]
holds for every $q\geq 1$ and $m\geq 1$.
Hence every element of 
\[
H^{0}(X_{m},{\cal O}_{X_{m}}(K_{X_{m}}+p_{m}^{*}(mK_{X}+L_{1} - K_{X}))
\otimes {\cal I}(p_{m}^{*}h^{m})\otimes {\cal O}_{X_{m}}/p_{m}^{*}{\cal I}_{X_{1}}) 
\]
extends to an element of 
\[
H^{0}(X_{m},{\cal O}_{X_{m}}(K_{X_{m}}+p_{m}^{*}(mK_{X}+L_{1} - K_{X}))
\otimes {\cal I}(p_{m}^{*}h^{m})).
\]
Also there exists a natural map 
\[
H^{0}(X_{1},{\cal O}_{X_{1}}(mK_{X}+L_{1})\otimes{\cal I}(h^{m}))
\rightarrow
\]
\[ 
\hspace{20mm} H^{0}(X_{m},{\cal O}_{X_{m}}(K_{X_{m}}+p_{m}^{*}(mK_{X}+L_{1} - K_{X}))
\otimes {\cal I}(p_{m}^{*}h^{m})\otimes {\cal O}_{X_{m}}/p_{m}^{*}{\cal I}_{X_{1}}). 
\]
Hence we can extend every element of
\[
p_{m}^{*}H^{0}(X_{1},{\cal O}_{X_{1}}(mK_{X}+L_{1})\otimes{\cal I}(h^{m}))
\]
to an element of 
\[
H^{0}(X_{m},{\cal O}_{X_{m}}(K_{X_{m}}+p_{m}^{*}(mK_{X}+L_{1} - K_{X}))
\otimes {\cal I}(p_{m}^{*}h^{m})).
\]
Since 
\[
H^{0}(X_{m},{\cal O}_{X_{m}}(K_{X_{m}}+p_{m}^{*}(mK_{X}+L_{1} - K_{X}))
\otimes {\cal I}(p_{m}^{*}h^{m}))\simeq 
\]
\[
H^{0}(X,{\cal O}_{X}(mK_{X}+L_{1})\otimes{\cal I}(h^{m}))
\]
holds by the isomorphism 
\[
p_{m*}({\cal O}_{X_{m}}(K_{X_{m}})\otimes {\cal I}(p_{m}^{*}h^{m}))
= {\cal O}(K_{X})\otimes {\cal I}(h^{m}),
\]
this completes the proof of Lemma 3.3.
\vspace{10mm}{\bf Q.E.D.} \\ 

Let $\tau$ be a general section in 
$H^{0}(X,{\cal O}_{X}(L_{1}))$.

Let $m_{1}$ be a sufficiently large positive integer 
and let $\sigma_{1}^{\prime}$
be a general element of 
\[
H^{0}(X_{1},{\cal O}_{X_{1}}(m_{1}K_{X})\otimes{\cal I}(h^{m_{1}})\otimes
{\cal M}_{x_{1},x_{2}}^{\lceil\sqrt[n_{1}]{\mu_{1}}(1-\varepsilon )\frac{m_{1}}
{\sqrt[n_{1}]{2}}\rceil}),
\]
where $x_{1},x_{2}\in X_{1}$ are distinct nonsingular points on $X_{1}$.

By Lemma 3.2, we may assume that $\sigma_{1}^{\prime}$ is nonzero.
Then by Lemma 3.3 we see that   
\[
\sigma_{1}^{\prime}\otimes\tau\in
H^{0}(X_{1},{\cal O}_{X_{1}}(m_{1}K_{X}+L_{1})\otimes{\cal I}(h^{m_{1}})\otimes
{\cal M}_{x_{1},x_{2}}^{\lceil\sqrt[n_{1}]{\mu_{1}}(1-\varepsilon )\frac{m_{1}}
{\sqrt[n_{1}]{2}}\rceil})
\]
extends to a section
\[
\sigma_{1}\in H^{0}(X,{\cal O}_{X}((m_{1}+\ell_{1} )K_{X})
\otimes{\cal I}(h^{m+\ell_{1}}))
\]
We may assume that  there exists a neighbourhood $U_{x,x^{\prime}}$ of $\{ x,x^{\prime}\}$ such that the divisor $(\sigma _{1})$  is smooth
on  $U_{x,x^{\prime}} - X_{1}$ by Bertini's theorem, if we take $\ell_{1}$
sufficiently large, since by Theorem 2.1, 
\[
H^{0}(X,{\cal O}_{X}(mK_{X}+L_{1})\otimes{\cal I}(h^{m}))
\rightarrow
H^{0}(X,{\cal O}_{X}(mK_{X}+L_{1})\otimes{\cal I}(h^{m}))/
{\cal O}_{X}(-X_{1})\cdot{\cal M}_{y})
\]
is surjective for every $y\in X$ and
 $m\geq 0$, where ${\cal O}_{X}(-X_{1})$
is the ideal sheaf of $X_{1}$.
We define a singular hermitian metric $h_{1}$ on $K_{X}$ by
\[
h_{1} = \frac{1}{\mid\sigma_{1}\mid^{\frac{2}{m_{1}+\ell_{1}}}}.
\]
Let $\varepsilon_{0}$ be a sufficiently small positive number and 
let ${\cal I}_{1}(\alpha )$ be the multiplier ideal sheaf of 
$h_{0}^{\alpha_{0}-\varepsilon_{0}}\cdot h_{1}^{\alpha}$,i.e.,
\[
{\cal I}_{1}(\alpha ) = {\cal L}^{2}({\cal O}_{X},
h_{0}^{\alpha_{0}-\varepsilon_{0}}h_{1}^{\alpha}/
h_{X}^{(\alpha_{0}+\alpha-\varepsilon_{0})}).
\]
Suppose that $x,x^{\prime}$ are nonsingular points on $X_{1}$.
Then we set $x_{1} = x, x_{2} = x^{\prime}$ and define $\alpha_{1}(=\alpha_{1}(x,y))> 0$ by
\[
\alpha_{1} := \inf\{\alpha\mid ({\cal O}_{X}/{\cal I}_{1}(\alpha ))_{x}
\neq 0\,\mbox{and}\, ({\cal O}_{X}/{\cal I}_{1}(\alpha ))_{x^{\prime}}\neq 0\}.
\]
By Lemma 3.3 we may assume that we have taken $m_{1}$ so that  
\[
\frac{\ell_{1}}{m_{1}} \leq 
\varepsilon_{0}\frac{\sqrt[n_{1}]{\mu_{1}}}{n_{1}\sqrt[n_{1}]{2}}
\]
holds.
\begin{lemma}
\[
\alpha_{1}\leq \frac{n_{1}\sqrt[n_{1}]{2}}{\sqrt[n_{1}]{\mu_{1}}} 
+ O(\varepsilon _{0})
\]
holds.
\end{lemma}
To prove Lemma 3.4, we need the following elementary lemma.
\begin{lemma}(\cite[p.12, Lemma 6]{t})
Let $a,b$ be  positive numbers. Then
\[
\int_{0}^{1}\frac{r_{2}^{2n_{1}-1}}{(r_{1}^{2}+r_{2}^{2a})^{b}}
dr_{2}
=
r_{1}^{\frac{2n_{1}}{a}-2b}\int_{0}^{r_{1}^{-{2}{a}}}
\frac{r_{3}^{2n_{1}-1}}{(1 + r_{3}^{2a})^{b}}dr_{3}
\]
holds, where 
\[
r_{3} = r_{2}/r_{1}^{1/a}.
\]
\end{lemma}
{\bf Proof of Lemma 3.3.}
Let $(z_{1},\ldots ,z_{n})$ be a local coordinate on a 
neighbourhood $U$ of $x$ in $X$ such that 
\[
U \cap X_{1} = 
\{ q\in U\mid z_{n_{1}+1}(q) =\cdots = z_{n}(q)=0\} .
\] 
We set $r_{1} = (\sum_{i=n_{1}+1}^{n}\mid z_{1}\mid^{2})^{1/2}$ and 
$r_{2} = (\sum_{i=1}^{n_{1}}\mid z_{i}\mid^{2})^{1/2}$.
Then there exists a positive constant $C$ such that 
\[
\parallel\sigma_{1}\parallel^{2}\leq 
C(r_{1}^{2}+r_{2}^{2\lceil\sqrt[n_{1}]{\mu_{1}}(1-\varepsilon )\frac{m_{1}}
{\sqrt[n_{1}]{2}}\rceil})
\]
holds on a neighbourhood of $x$, 
where $\parallel\,\,\,\,\parallel$ denotes the norm with 
respect to $h_{X}^{m_{1}+\ell_{1}}$.
We note that there exists a positive integer $M$ such that 
\[
\parallel\sigma\parallel^{-2} = O(r_{1}^{-M})
\]
holds on a neighbourhood of the generic point of $U\cap X_{1}$,
where $\parallel\,\,\,\,\parallel$ denotes the norm with respect to 
$h_{X}^{m_{0}}$. 
Then by Lemma 3.5, we have the inequality 
\[
\alpha_{1} \leq (\frac{m_{1}+\ell_{1}}{m_{1}})\frac{n_{1}\sqrt[n_{1}]{2}}{\sqrt[n_{1}]{\mu_{1}}} 
+ O(\varepsilon _{0})
\] 
holds. 
By using the fact that 
\[
\frac{\ell_{1}}{m_{1}} \leq 
\varepsilon_{0}\frac{\sqrt[n_{1}]{\mu_{1}}}{n_{1}\sqrt[n_{1}]{2}}
\]
we obtain that 
\[
\alpha_{1}\leq \frac{n_{1}\sqrt[n_{1}]{2}}{\sqrt[n_{1}]{\mu_{1}}} 
+ O(\varepsilon _{0})
\]
holds.
{\bf Q.E.D.} \vspace{10mm} \\
If $x$ or $x^{\prime}$ is a singular point on $X_{1}$, we need the following lemma.
\begin{lemma}
Let $\varphi$ be a plurisubharmonic function on $\Delta^{n}\times{\Delta}$.
Let $\varphi_{t}(t\in\Delta )$ be the restriction of $\varphi$ on
$\Delta^{n}\times\{ t\}$.
Assume that $e^{-\varphi_{t}}$ does not belong to $L^{1}_{loc}(\Delta^{n},O)$
for every $t\in \Delta^{*}$.

Then $e^{-\varphi_{0}}$ is not locally integrable at $O\in\Delta^{n}$.
\end{lemma}
Lemma 3.6 is an immediate consequence of \cite{o-t}.
Using Lemma 3.6 and Lemma 3.5, we see that Lemma 3.4 holds
by letting $x_{1}\rightarrow x$ and $x_{2}\rightarrow x^{\prime}$.

\vspace{10mm}

For the next we consider Case 1.2 and Case 2.2.  
We note that in Case 2.2 by modifying $\sigma$ a little bit 
, if necessary we may assume that
$({\cal O}_{X}/{\cal I}(\alpha_{0}-\varepsilon ))_{x^{\prime}}\neq 0$ 
and $({\cal O}_{X}/{\cal I}(\alpha_{0}-\varepsilon^{\prime}))_{x} = 0$ hold
for a sufficiently small positive number $\varepsilon^{\prime}$. 
For example it is sufficient to replace $\sigma$ by 
the following $\sigma^{\prime}$ constructed below.

Let $X^{\prime}_{1}$ be a subvariety which defines a branch of 
\[
\mbox{Spec}({\cal O}_{X}/{\cal I}(\alpha +\delta))
\]
at $x^{\prime}$.  By the assumption (changing $X_{1}$, if necessary) we may assume that $X_{1}^{\prime}$ does not 
contain $x$.  
Let $m^{\prime}$ be a sufficiently large positive integer such that 
$m^{\prime}/m_{0}$ is sufficiently small (we can take $m_{0}$ 
arbitrary large). 

Let $\tau_{x^{\prime}}$ be a general element of 
\[
H^{0}(X,{\cal O}_{X}(m^{\prime}K_{X})\otimes 
{\cal I}_{X_{1}^{\prime}}), 
\]
where ${\cal I}_{X_{1}^{\prime}}$ is the ideal sheaf of 
$X_{1}^{\prime}$. 
If we take $m^{\prime}$ sufficiently large,
 $\tau_{x^{\prime}}$ is not identically zero. 
We set 
\[
\sigma^{\prime} = \sigma\cdot\tau_{x^{\prime}}. 
\] 
Then we see that the new singular hermitian metric $h_{0}^{\prime}$
defined by $\sigma^{\prime}$ satisfies the desired property.

In these cases, instead of Lemma 3.2, we use the following simpler lemma.

\begin{lemma} Let $\varepsilon$ be a sufficiently small positive number and let $x_{1}$ be a smooth point on $X_{1}$. 
Then for a sufficiently large $m >1$,
\[
H^{0}(X_{1},{\cal O}_{X_{1}}(mK_{X})\otimes{\cal I}(h^{m})\otimes
{\cal M}_{x_{1}}^{\lceil\sqrt[n_{1}]{\mu_{1}}(1-\varepsilon )m\rceil})\neq 0
\]
holds.
\end{lemma}

Then taking a general $\sigma_{1}^{\prime}$ in
\[
H^{0}(X_{1},{\cal O}_{X_{1}}(m_{1}K_{X})\otimes{\cal I}(h^{m_{1}})\otimes
{\cal M}_{x_{1}}^{\lceil\sqrt[n_{1}]{\mu_{1}}(1-\varepsilon )m_{1}
\rceil}),
\]
for a sufficiently large $m_{1}$.
As in Case 1.1 and Case 2.1 we obtain a proper subvariety
$X_{2}$ in $X_{1}$ also in this case.

Inductively for distinct points $x,x^{\prime}\in X^{\circ}$, we construct a strictly decreasing
sequence of subvarieties
\[
X = X_{0}\supset X_{1}\supset \cdots \supset X_{r}\supset X_{r+1} = \{ x\}\cup 
R_{x^{\prime}}\,\mbox{or}\, \{ x^{\prime}\}\cup R_{x},
\]
where $R_{x^{\prime}}$ (or $R_{x}$) is a subvariety such that $x$ deos not
belong to $R_{x^{\prime}}$ and $x^{\prime}$ belongs to $R_{x^{\prime}}$.
and invariants (depending on small positive numbers $\varepsilon_{0},\ldots ,
\varepsilon_{r-1}$, large positive integers $m_{0},m_{1},\ldots ,m_{r}$, etc.) :
\[
\alpha_{0} ,\alpha_{1},\ldots ,\alpha_{r},
\]
\[
\mu_{0},\mu_{1},\ldots ,\mu_{r}
\]
and
\[
n >  n_{1}> \cdots > n_{r}.
\]
By Nadel's vanishing theorem we have the following lemma.
\begin{lemma} 
Let $x,x^{\prime}$ be two distinct points on $X^{\circ}$. 
Then for every $m\geq \lceil\sum_{i=0}^{r}\alpha_{i}\rceil +1$,
$\Phi_{\mid mK_{X}\mid}$ separates $x$ and $x^{\prime}$.
\end{lemma}
{\bf Proof}. 
Let us define the singular hermitian metric $h_{x,x^{\prime}}$ of $(m-1)K_{X}$ defined by  
\[
h_{x,x^{\prime}} = (\prod_{i=0}^{r-1}h_{i}^{\alpha_{i}-\varepsilon_{i}})\cdot
 h_{r}^{\alpha_{r}+\varepsilon_{r}}\cdot h^{(m-1-(\sum_{i=0}^{r-1}(\alpha_{i}-\varepsilon_{i}))- (\alpha_{r}+\varepsilon_{r})-\delta_{L})}\cdot h_{L}^{\delta_{L}},
\]
where $h_{L}$ is a $C^{\infty}$-hermitian metric on the {\bf Q}-line bundle 
$L := K_{X}-E$ with strictly positive curvature and $\delta_{L}$ be a sufficiently small positive number. 
Then we see that  ${\cal I}(h_{x,x^{\prime}})$ defines a subscheme of 
$X$ with isolated support around $x$ or $x^{\prime}$ by the definition of 
the invariants $\{\alpha_{i}\}$'s. 
By the construction the curvature current $\Theta_{h_{x,x^{\prime}}}$ is strictly positive on $X$. 
Then by Nadel's vanishing theorem (Theorem 2.1) we see that 
\[
H^{1}(X,{\cal O}_{X}(mK_{X})\otimes {\cal I}(h_{x,x^{\prime}})) = 0.
\]
This implies that $\Phi_{\mid mK_{X}\mid}$ separates 
$x$ and $x^{\prime}$.   {\bf Q.E.D.} 

\subsection{Construction of the stratification as a family}

In this subsection we shall construct the above stratification as a family. 

We note that for a fixed pair $(x,x^{\prime}) \in X^{\circ}\times X^{\circ}-\Delta_{X}$, $\sum_{i=0}^{r}\alpha_{i}$ depends on the choice of $\{ X_{i}\}$'s, where 
$\Delta_{X}$ denotes the diagonal of $X\times X$. 
Moving $(x,x^{\prime})$ in $X^{\circ}\times X^{\circ} - \Delta_{X}$, 
we  shall consider the above operation simultaneously.
Let us explain the procedure. 
We set 
\[
B := X^{\circ}\times X^{\circ} - \Delta_{X}.
\] 
Let 
\[
p : X\times B\longrightarrow X
\]
be the first projection and let   
\[
q : X\times {B}
\longrightarrow B
\]
be the second projection. 
Let $Z$ be the subvariety of $X\times B$ defined by
\[
Z := \{ (x_{1},x_{2},x_{3}) : X\times B \mid 
x_{1} = x_{2} \,\,\mbox{or} \,\, x_{1} = x_{3} \} .
\]
In this case we consider 
\[
q_{*}{\cal O}_{X\times B}(m_{0}p^{*}K_{X})
\otimes {\cal I}_{Z}^{\lceil\sqrt[n]{\mu_{0}}(1-\varepsilon )\frac{m_{0}}{\sqrt[n]{2}}\rceil }
\]
instead of 
\[
H^{0}(X,{\cal O}_{X}(m_{0}K_{X})\otimes
{\cal M}_{x,x^{\prime}}^{\lceil\sqrt[n]{\mu_{0}}(1-\varepsilon )\frac{m_{0}}{\sqrt[n]{2}}\rceil}),
\]
where ${\cal I}_{Z}$ denotes the ideal sheaf of $Z$. 
Let $\tilde{\sigma}_{0}$ be a nonzero global meromorphic 
section  of 
\[
q_{*}{\cal O}_{X\times B}(m_{0}p^{*}K_{X})
\otimes {\cal I}_{Z}^{\lceil\sqrt[n]{\mu_{0}}(1-\varepsilon )\frac{m_{0}}{\sqrt[n]{2}}\rceil } 
\]
on $B$ for a sufficiently large positive integer $m_{0}$.
We define the singular hermitian metric $\tilde{h}_{0}$ 
on $p^{*}K_{X}$ by 
\[
\tilde{h}_{0}:= \frac{1}{\mid \tilde{\sigma}_{0}\mid^{2/m_{0}}}.
\]
We shall replace $\alpha_{0}$  by 
\[
\tilde{\alpha}_{0} 
:= \inf \{\alpha > 0\mid 
\mbox{the generic point of}\,\, Z \subseteq \mbox{Spec}
({\cal O}_{X \times B}/{\cal I}(h_{0}^{\alpha}))\} .
\]
Then for every $0 < \delta << 1$, there exists a Zariski 
open subset $U$ of $B$ such that for every $b \in U$, 
$\tilde{h}_{0}\mid_{X\times\{ b\}}$ is well defined and  
\[
b \not{\subseteq}\mbox{Spec}({\cal O}_{X\times\{ b\}}/{\cal I}(\tilde{h}_{0}^{\alpha_{0}-\delta}\mid_{X\times\{ b\}})),
\]
where we have identified $b$ with distinct two points in $X$. 
And also by Lemma 3.6, we see that 
\[
b \subseteq \mbox{Spec}({\cal O}_{X\times\{ b\}}/{\cal I}(\tilde{h}_{0}^{\alpha_{0}}\mid_{X\times\{ b\}})),
\]
holds for every $b\in B$. 
Let $\tilde{X}_{1}$ be an irreducible component of  
\[
\mbox{Spec}({\cal O}_{X\times B}/{\cal I}(\tilde{h}_{0}^{\alpha_{0}}))
\]
containing $Z$. 
We note that $\tilde{X}_{1}\cap q^{-1}(b)$ may not be 
irreducible even for a general $b\in B$. 
But if we take a suitable finite cover 
\[
\phi_{0} : B_{0} \longrightarrow B,
\]
on the base change $X\times_{B}B_{0}$, $\tilde{X}_{1}$ 
defines a family of irreducible subvarieties
\[
f_{1} : \hat{X}_{1} \longrightarrow U_{0}
\]
of $X$ parametrized by a nonempty Zariski open subset
 $U_{0}$ of $\phi_{0}^{-1}(U)$.
We set 
\[
\tilde{\mu}_{1} := \inf_{b_{0}\in U_{0}}\mu (f_{1}^{-1}(b_{0}),(K_{X},h)).
\]
We note that by its definition the volume $\mu (f_{1}^{-1}(b_{0}),(K_{X},h))$ 
is constant on a nonempty open subset say $U_{0}^{\prime}$ of $U_{0}$  
with respect to countable Zariski topology.
We denote the constant by $\tilde{\mu}_{0}$. 
Continueing this process 
we may construct a finite morphism 
\[
\phi_{r} : B_{r} \longrightarrow B
\]
and a nonempty Zariski open subset $U_{r}$ of $B_{r}$ 
which parametrizes a family of stratification 
\[
X \supset  X_{1} \supset X_{2} \supset \cdots \supset X_{r} 
\supset X_{r+1} = \{ x\} \cup R_{x^{\prime}} (\mbox{resp.} \,\,
\{ x^{\prime}\} \cup R_{x})  
\] 
constructed as before, 
where $R_{x}$ (resp. $R_{x^{\prime}}$) is a subvariety of $X$ which is 
disjoint from $x^{\prime}$ (resp. $x$). 
And we also obtain invariants $\{\tilde{\alpha}_{0},
\ldots ,\tilde{\alpha}_{r}\}$, $\{\tilde{\mu}_{0},\ldots ,\tilde{\mu}_{r}\}$,
$\{ n = \tilde{n}_{0}\ldots ,\tilde{n}_{r}\}$.
Hereafter we denote these invariants without $\,\,\tilde{} \,\,$ for simplicity.  By the same proof as Lemma 3.4, we have the following lemma. 
\begin{lemma}
\[
\alpha_{i}\leq \frac{n_{i}\sqrt[n_{i}]{2}}{\sqrt[n_{i}]{\mu_{i}}} + O(\varepsilon_{i-1})
\]
hold for $1\leq i\leq r$.
\end{lemma}
By Lemma 3.8 we obtain the following proposition. 
\begin{proposition}
For every 
\[
m > \lceil\sum_{i=0}^{r}\alpha_{i}\rceil + 1
\]
$\mid mK_{X}\mid$ gives a birational rational map from $X$ into 
a projective space.
\end{proposition}

\begin{lemma} Let $X_{i}$ be a strata of a very general member of 
the stratification parametrized by $B_{r}$. 
If $\Phi_{\mid mK_{X}\mid}\mid_{X_{i}}$ is birational rational map
onto its image, then
\[
\deg \Phi_{\mid mK_{X}\mid}(X_{i})\leq m^{n_{i}}\mu_{i}
\]
holds.
\end{lemma}
{\bf Proof}.
Let $p : \tilde{X}\longrightarrow X$ be the resolution of 
the base locus of $\mid mK_{X}\mid$ and let 
\[
p^{*}\mid mK_{X}\mid = \mid P_{m}\mid + F_{m}
\]
be the decomposition into the free part $\mid P_{m}\mid$ 
and the fixed component $F_{m}$. 
Let $p_{i} : \tilde{X}_{i}\longrightarrow X_{i}$ be the resolution
of the base locus of $\Phi_{\mid mK_{X}\mid}\mid_{X_{i}}$ 
obtained by the restriction of $p$ on $p^{-1}(X_{i})$. 
Let 
\[
p_{i}^{*}(\mid mK_{X}\mid_{X_{i}}) = \mid P_{m,i}\mid + F_{m,i}
\]
be the decomposition into the free part $\mid P_{m,i}\mid$ and 
the fixed part $F_{m,i}$.
We have
\[
\deg \Phi_{\mid mK_{X}\mid}(X_{i}) = P_{m,i}^{n_{i}}
\]
holds.
Then by the ring structure of $R(X,K_{X})$, we have an injection 
\[
H^{0}(\tilde{X},{\cal O}_{\tilde{X}}(\nu P_{m}))\rightarrow 
H^{0}(X,{\cal O}_{X}(m\nu K_{X})\otimes{\cal I}(h^{m\nu}))
\]
for every $\nu\geq 1$, since 
the righthandside is isomorphic to 
$H^{0}(X,{\cal O}_{X}(m\nu K_{X}))$ by the definition of 
an AZD.
We note that since ${\cal O}_{\tilde{X}}(\nu P_{m})$ is globally generated
on $\tilde{X}$, for every $\nu \geq 1$ we have the injection 
\[
{\cal O}_{\tilde{X}}(\nu P_{m})\rightarrow p^{*}({\cal O}_{X}(m\nu K_{X})\otimes{\cal I}(h^{m\nu})).
\]
Hence there exists a natural morphism 
\[
H^{0}(\tilde{X}_{i},{\cal O}_{\tilde{X}_{i}}(\nu P_{m,i}))
\rightarrow 
H^{0}(X_{i},{\cal O}_{X_{i}}(m\nu K_{X})\otimes{\cal I}(h^{m\nu})/\mbox{tor})
\]
for every $\nu\geq 1$. 
This morphism is clearly injective. 
This implies that 
\[
\mu_{i} \geq  m^{-n_{i}}\mu (\tilde{X}_{i},P_{m,i})
\]
holds. 
Since $P_{m,i}$ is nef and big on $X_{i}$ we see that 
\[
\mu (\tilde{X}_{i},P_{m,i}) = P_{m,i}^{n_{i}}
\]
holds.
Hence
\[
\mu_{i}\geq m^{-n_{i}}P_{m,i}^{n_{i}}
\]
holds.  This implies that
\[
\deg \Phi_{\mid mK_{X}\mid}(X_{i})\leq \mu_{i}\cdot m^{n_{i}}
\]
holds.
{\bf Q.E.D.}

\section{Proof of Theorem 1.2}
Let
\[
X \supset  X_{1} \supset X_{2} \supset \cdots \supset X_{r} 
\supset X_{r+1} = \{ x\} \cup R_{x^{\prime}} (\mbox{resp.} \,\,
\{ x^{\prime}\} \cup R_{x})  
\] 
be a very general stratification constructed as 
in the last section. 
 
Suppose that $\mid \lceil (\lceil\sum_{\ell=1}^{n}\sqrt[\ell]{2}\ell\rceil +1)K_{X}\mid$ does not give a birational
rational map from $X$ into a projective space.  
Then  by Proposition 3.1, we see that
\[
\max_{i}\frac{\alpha_{i}}{\sqrt[n_{i}]{2}n_{i}} \geq 1
\]
holds.
Let $k$ be the number such that 
\[
\frac{\alpha_{k}}{\sqrt[n_{k}]{2}n_{k}} = \max_{i}\frac{\alpha_{i}}{\sqrt[n_{i}]{2}n_{i}}.
\]
By Lemma 3.9 we see that 
\[
\mu_{k} < 1
\]
holds. 
We set 
\[
\alpha := \lceil\sum_{i=0}^{r} \alpha_{i}\rceil + 1.
\]
Now we see that by Lemma 3.10 and Lemma 3.9 
\[
\mbox{deg}\,\Phi_{\mid \alpha K_{X}\mid}(X_{k})
\leq (\lceil (\sum_{\ell=1}^{n}\frac{\sqrt[\ell]{2}\,\,\ell}{\sqrt[n_{k}]{2}\,\,n_{k}})\alpha_{k}\rceil +1)^{n_{k}}\mu_{k}
\]
\[
\,\,\,\,\,\,\leq ((\lceil\sum_{\ell=1}^{n}\sqrt[\ell]{2}\,\,\ell\rceil +1)
\frac{1}{\sqrt[n_{k}]{\mu_{k}}})^{n_{k}}\mu_{k}
\]
\[
\,\,\,\,\,\,\leq (\lceil\sum_{\ell=1}^{n}\sqrt[\ell]{2}\,\,\ell\rceil +1)^{n_{k}}
\]
hold.
Since such $\{ X_{k}\}$ form a dominant family of 
subvarieties on $X$ by the construction of 
the stratifications, this completes the proof of Theorem 1.2.

\section{Proof of Theorem 1.3}
Let $X$ be a smooth projective $n$-fold of general type and suppose that 
$\mid (\lceil\sum_{\ell=1}^{n}\sqrt[\ell]{2}\,\,\ell\rceil +1)K_{X}\mid$ 
does not give a birational embedding. 
Then by Theorem 1.2 and its proof, there exists some $1\leq d\leq n$ such that 
there exists a dominant family of subvarieties 
\[
\varpi : {\cal X}_{k}\longrightarrow S_{k}
\]
which parametrizes the strata $X_{k}$ of a general stratification 
\[
X \supset  X_{1} \supset X_{2} \supset \cdots \supset X_{r} 
\supset X_{r+1} = \{ x\} \cup R_{x^{\prime}} (\mbox{resp.} \,\,
\{ x^{\prime}\} \cup R_{x})
\] 
such that for 
$\alpha := \lceil\sum_{i=0}^{r} \alpha_{i}\rceil + 1$
\[
\deg \Phi_{\mid\alpha K_{X}\mid}(X_{k}) \leq 
(\lceil\sum_{\ell=1}^{n}\sqrt[\ell]{2}\,\,\ell\rceil +1)^{d^{\prime}}
\]
holds where $d^{\prime} = \dim X_{k}$. 
Let 
\[
p : {\cal X}_{k}\longrightarrow X
\]
be the natural morphism.
Inductively we define a sequence of (possibly reducible) subvarieties $F_{i} (i\geq 0)$ by 
\[
F_{0} = X_{k} ([X_{k}]\in S_{k})
\]
and for $i\geq 0$
\[
F_{i+1} = \mbox{the closure of}\,\,\,\,
p(\varpi^{-1}(\varpi (\pi^{-1}(\mbox{the generic points of $F_{i}$})))).
\]
Then $\{ F_{i}\}_{i\geq 0}$ is increasing and by the Noetherian property, we see that there exists some $\ell \geq 0$ such that 
\[
\dim F_{\ell} = \dim F_{\ell^{\prime}}
\]
for every $\ell^{\prime}\geq \ell$.
Let $F_{\ell,0}$ be a maximal dimensional component of $F_{\ell}$. 

If we start from a general $[X_{k}]\in S_{k}$ and choose $F_{\ell,0}$ properly,
we may assume that 
$\{ F_{\ell,0}\}$ form a family. 
We note that possibly $F_{\ell,0}$ may not be determined 
only by $X_{k}$ because of a monodoromy phenomenon. 
Let 
\[
\varpi_{0} : {\cal U}_{0} \longrightarrow T_{0}
\]
be the family of such $\{ F_{\ell,0}\}$.
Then it is again dominant. 
Let 
\[
p_{0}: {\cal U}_{0}\longrightarrow X
\]
be the natural morphism. 
We see that $p_{0}$ is birational, since for a general $x\in X$, 
$p_{0}^{-1}(x)$ is a point (otherwise it contradicts to the maximality of 
$\dim F_{\ell,0}$). 
Hence $\varpi_{0}$ induces a rational fibration structure
\[
f : X -\cdots\rightarrow Y.
\] 
Let $F$ be a general fiber of $f$.
To completes the proof of Theorem 1.3 we need the following lemma.  
\begin{lemma}
\[
\deg \Phi_{\mid\alpha K_{X}\mid}(F) 
\leq 
(\lceil\sum_{\ell=1}^{n}\sqrt[\ell]{2}\,\,\ell\rceil +1)^{d^{2}}d^{d}
\]
\end{lemma}
{\em Proof}.  
By taking a suitable modification of $X$, we may assume that 
the following conditions are satisfied :
\begin{enumerate}
\item $\Phi_{\mid\alpha K_{X}\mid}$ is a morphism on $X$,
\item there exists a regular fibration 
\[
f : X\longrightarrow Y
\]
induced by $p_{0}: {\cal U}_{0}\longrightarrow X$ as above. 
\end{enumerate}
Let $\mid H\mid$ be the free part of $\mid \alpha K_{X}\mid$.
We set 
\[
a := (\lceil\sum_{\ell=1}^{n}\sqrt[\ell]{2}\,\,\ell\rceil +1)^{d^{\prime}}.
\]
Let $F$ be a general fiber of $f$.
Suppose that 
\[
\deg \Phi_{\mid\alpha K_{X}\mid}(F) > a^{d}d^{d}
\]
holds.  
\begin{lemma}
Let us fix an arbitrary point $x_{0}$ on $F$.  
Then for a sufficiently large $m$ there exists a section
\[
\sigma \in \Gamma (F,{\cal O}_{F}(mH))-\{ 0\}
\]
such that 
\[
\mbox{mult}_{x_{0}}(\sigma )  > mad +1
\]
holds.
\end{lemma}
{\bf Proof}. 
Since 
\[
\dim {\cal O}_{F}/{\cal M}_{x_{0}}^{m} 
= \left(\begin{array}{c} d+m-1 \\
d \end{array}\right) = \frac{1}{d!}m^{d} + O(m^{d-1}) 
\]
and 
\[
\dim H^{0}(F,{\cal O}_{X}(mH)) 
= \frac{1}{d!}a^{d}d^{d}m^{d} + O(m^{d-1})
\]
hold, the lemma is clear. 
{\bf Q.E.D.} \vspace{5mm} \\ 

Let 
\[
\varpi_{F} : {\cal X}_{k}(F) \longrightarrow S_{k}(F)
\]
be the family of the strata $X_{k}$ contained in $F$.
Then there exists a subvariety $S^{\prime}_{k}(F)$ such that 
\begin{enumerate}
\item $\dim S^{\prime}_{k}(F) = \dim F - \dim X_{k}$,
\item  We set $\tilde{F}: = \varpi_{F}^{-1}(S^{\prime}_{k}(F))$.
Then $p_{\tilde{F}} : \tilde{F}\longrightarrow F$ 
is generically finite. 
\end{enumerate}
Then 
\[
\tilde{\varpi}_{F} : \tilde{F} \longrightarrow S^{\prime}_{k}(F)
\]
is an algebraic fiber space. 
We set 
\[
W_{j} = \{ \tilde{x}\in \tilde{F}\mid \mbox{mult}_{\tilde{x}}\,\,p_{\tilde{F}}^{*}(\sigma )
 \geq 1+maj\} .
\]
We have a decending chain of subvarieties
\[
W_{0}\supset W_{1}\supset \cdots \supset W_{d}\ni p_{\tilde{F}}^{-1}(x_{0}).
\]
For each $j\in \{ 0,\ldots ,n\}$, we choose an irreducible component of $W_{j}$ containing $x$, and denote this irreducible component by $W_{j}^{\prime}$.
We may assume that these irreducible components have been chosen so that 
\[
W_{0}^{\prime}\supset W_{1}^{\prime}\supset \cdots \supset W_{d}^{\prime}\ni  p_{\tilde{F}}^{-1}(x_{0}).
\]
Since this chain has length greater than $d = \dim F$, there exists some 
$j$ such that 
\[
W_{j}^{\prime} = W_{j+1}^{\prime}
\]
holds. 
We set $W = W_{j}^{\prime}$. 

Let ${\cal C}$ denote the family of irreducible curves
\[
C = H_{1}\cap\cdots \cap H_{d^{\prime}-1}\cap X_{k}^{\prime}
\] 
on $F$ which is obtained as the  intersection of 
$(d^{\prime}-1)$-members $H_{1},\ldots H_{d^{\prime}-1}$ of $\mid H\mid$  
and a strata $X_{k}^{\prime}([X_{k}^{\prime}] \in S_{k})$ 
which is contained in $F$. 
We note that for a general member $C$ of ${\cal C}$, 
the inverse image of $C$ in the normalization of $X_{k}^{\prime}$ 
is smooth.  
Hence a general member of ${\cal C}$ is immersed in $X$ that means the 
differential of the natural morphism from the normalization of $C$ to
$X$ is nowhere vanishing.  
By the construction ${\cal C}$ determines 
a dominant family of 
curves $\tilde{{\cal C}}$ on $\tilde{F}$.
We note that by the construction of $F$
 a member $X_{k}^{\prime}$ of $S_{k}$ intersects 
$F$, then $X_{k}^{\prime}$ should be contained 
in $F$.
Now we note have that 
\[
H\cdot C \leq  a
\] 
holds.

Now we quote the following two theorems (the satements are slightly generalized, 
but the proofs are completely same). 
\begin{theorem}(\cite[p. 686,Theorem 2]{n2})
Fix a positve integers $\ell$ and $k$.
Let $f : M\longrightarrow \Delta^{k}$ be a smooth family 
of irreducible curves. 
Let $L$ be a holomorphic line bundle on $M$ such that the restriction
of $L$ to $f^{-1}(0)$ has degree $\ell$ and let $s\in H^{0}(M,{\cal O}_{M}(L))$.
Let $V_{j}(s)$ denote the complex subspace of $M$ consinsting precisely those points at which the vanishing order of $s$ is at least $j$. 

Then either
\[
f^{-1}(0) \cap V_{j+\ell}(s) = \emptyset
\]
or 
\[
f^{-1}(0) \subset V_{j}(s)
\]
holds. 
\end{theorem}

\begin{theorem}(\cite[p.686,Theorem 3]{n2})
Let $M$ be a complex manifold, let $L$ be a holomorphic line bundle on $M$
and let $s\in H^{0}(M,{\cal O}_{M}(L))$. 
Let ${\cal N}$ be an irreducible family of immersed curves thich is free and set 
$\ell := \deg_{N}L$, where $N$ is a member of ${\cal N}$. 
For any integer $j$, either
\[
N \cap V_{j+\ell}(s) = \emptyset
\]
or
\[
N\subset V_{j}(s)
\]
holds.
\end{theorem}

\begin{lemma}
If $\tilde{C}$ be an immersed member of $\tilde{\cal C}$ which intersects 
$W_{j}$, then $\tilde{C}$ is contained in $W_{j-1}$.
\end{lemma} 
{\bf Proof}. 
We note that $\tilde{C}$ is free by the construction.
Then the lemma follows from Theorem 5.2.  
{\bf Q.E.D.} \vspace{5mm} \\

Now we fix a general $W$. 
Then by the definition of ${\cal C}$ and Lemma 5.1 implies 
that if  for a fiber $\tilde{X}_{k}$ of
\[
\tilde{\varpi}_{F} : \tilde{F} \longrightarrow S^{\prime}_{k}(F)  
\]
\[
\tilde{X}_{k} \cap W \neq \emptyset
\]
holds, then 
\[
\tilde{X}_{k}\subseteq W
\]
holds. 
We note that $W$ is defined by the pullback of the section of $mH$ on $F$. 
By the inductive construction of $F$, 
we see that if we take $W$ general,  $W = \tilde{F}$ holds.
This is the contradiction. 
Hence we see that 
\[
\deg \Phi_{\mid\alpha K_{X}\mid}(F) 
\leq a^{d}d^{d}
\]
holds.  
Since 
\[
a^{d}d^{d} \leq (\lceil\sum_{\ell=1}^{n}\sqrt[\ell]{2}\,\,\ell\rceil +1)^{d^{2}}d^{d}
\]
holds, this completes the proof of Theorem 1.3.
 
Author's address\\
Hajime Tsuji\\
Department of Mathematics\\
Tokyo Institute of Technology\\
2-12-1 Ohokayama, Megro 152\\
Japan \\
e-mail address: tsuji@math.titech.ac.jp


\begin{thebibliography}{99}
\bibitem{a-s} U. Anghern-Y.-T. Siu, Effective freeness and point separation
for adjoint bundles, Invent. Math. 122 (1995), 291-308.
\bibitem{b-d} T. Bandman and G. Dethloff, Estimates of the number of rational mappings from a fixed variety to varieties of general type, 
Ann. Institute Fourier 47 (1997), 801-824. 
\bibitem{b} E. Bombieri, Algebraic values of meromoprhic maps, Invent. Math. 10
(1970), 267-287.
\bibitem{b2} E. Bombieri, Addendum to my paper: Algebraic values of 
meromorphic maps, Invent. Math. 11, 163-166. 
\bibitem{b3} E. Bombieri, Canonical models of surfaces of general type,
Publ. I.H.E.S. 42 (1972), 171-219.
\bibitem{ca} Mark Andrea A. de Catalado, Singular hermitian metrics on vector bundles, math.AG/9708003, J. fur Reine Angewande Math. 502(1998), 93-102. 
\bibitem{he} S. Helmke, On Fujita's conjecture, Duke Math. J. 88(1997),
 201-216.
\bibitem{h} L. H\"{o}rmander, An Introduction to Complex Analysis in Several
Variables 3-rd ed.,North-Holland(1990).
\bibitem{l} P. Lelong, Founctions Plurisousharmoniques et Formes 
Differentielles Positives, Gordon and Breach (1968).
\bibitem{luo} T. Luo, Global 2-forms on regular threefolds of general type, 
Duke Math. J. 71 (1993), 859-869.
\bibitem{luo2} T. Luo, Plurigenera of regular threefolds, Math. Z. 217 (1994),
37-46.
\bibitem{n}A.M. Nadel, Multiplier ideal shaves and existence of K\"{a}hler-Einstein
metrics of positive scalar curvature, Ann. of Math. 132 (1990),549-596.
\bibitem{n2} A.M. Nadel, The boundedness of degree of Fano varieties with Picard number one, J. of Amer. Math. Soc. 4 (1991), 681-692. 
\bibitem{o-t}T. Ohsawa and K. Takegoshi, $L^{2}$-extention of holomorphic
functions, Math. Z. 195 (1987),197-204.
\bibitem{si} Y.-T. Siu, Analyticity of sets associated to Lelong numbers
and the extension of closed positive currents, Invent. Math. 27 (1974), 
53-156.
\bibitem{tu}H. Tsuji, Analytic Zariski decomposition, Proc. of Japan Acad.
61(1992) 161-163.
\bibitem{tu1} H. Tsuji, Existence and applications of Analytiv Zariski decompositions, Geometry and Analysis in Several Complex Variables, Trends in Math.  (1999), 253-271. 
\bibitem{tu2} H. Tsuji, On the structure of pluricanonical systems of projective varieties of general type, TIT preprint series (1997). 
\bibitem{tu3} H. Tsuji, Finite generation of canonical rings, math.AG/9908083 (1999).
\bibitem{tu4} H. Tsuji, Pluricanonical systems of projective varieties of 
general type, math.AG/9909021 (1999). 
\bibitem{t}H. Tsuji, Global generation of adjoint bundles, Nagoya Math. J.
142 (1996),5-16.
\end{thebibliography}
\end{document}